\documentclass{artjlt-Preprint}

\usepackage{amsmath,amsfonts,amssymb,dsfont}
\usepackage{hyperref,xfrac,braket}

\providecommand{\R}{\mathbb{R}}
\renewcommand{\H}{{\mathbb{H}}}
\providecommand{\C}{\mathbb{C}}

\providecommand{\D}{\mathbb{D}}

\providecommand{\tr}{\mathrm{tr}}

\renewcommand{\ker}{\mathrm{ker}}

\renewcommand{\H}{\mathbb{H}}

\newcommand{\SL}{\mathrm{SL}}
\newcommand{\SO}{\mathrm{SO}}

\newcommand{\GL}{\mathrm{GL}}

\newcommand{\1}{\mathds{1}}

\title{Algebraically Independent Generators for the Algebra of Invariant Differential Operators on $\mathrm{SL}_n(\mathbb R)/\mathrm{SO}_n(\mathbb R)$}                                     
\author{Dominik Brennecken, Lorenzo Ciardo and Joachim Hilgert} 
\lastname{Brennecken, Ciardo and Hilgert}  


\address{%
Dominik Brennecken\\               
Institut f\"ur Mathematik\\ 
Universit\"at Paderborn\\ 
D-33095 Paderborn\\ 
Germany\\ 
bdominik@math.upb.de 
}

\address{%
Lorenzo Ciardo\\               
Department of Mathematics\\ 
University of Oslo\\ 
Postboks 1053, Blindern 0316, Oslo\\ 
Norway\\ 
lorenzci@math.uio.no 
}

\address{%
Joachim Hilgert\\               
Institut f\"ur Mathematik\\ 
Universit\"at Paderborn\\ 
D-33095 Paderborn\\ 
Germany\\            
hilgert@upb.de                
}
\begin{document}


\maketitle

\begin{abstract}
 We provide an explicit set of algebraically independent generators for the algebra of invariant differential operators on the Riemannian symmetric space associated with $\SL_n(\R)$.
\end{abstract}


\section{Introduction}

Let $\mathrm{Pos}_n(\R)=\GL_n(\R)/\mathrm{O}_n(\R)$ be the space of positive definite real $n\times n$-matrices and $\mathrm{SPos}_n(\R)=\SL_n(\R)/\mathrm{SO}_n(\R)$ its subset of elements of determinant $1$. Both spaces are Riemannian symmetric spaces whose algebras $\mathbb D(\mathrm{Pos}_n(\R))$, respectively $\mathbb D(\mathrm{SPos}_n(\R))$, of invariant differential operators are commutative. The Harish-Chandra isomorphism together with Chevalley's Theorem shows that these algebras are isomorphic to the polynomial algebras in $n$, respectively $n-1$, variables. For  $\mathbb D(\mathrm{Pos}_n(\R))$ one finds various explicit sets of algebraically independent generators in the literature. One example of such a set is given by the Maass-Selberg operators $\delta_1,\ldots,\delta_n$ described in detail in \cite{JL05}. It is much harder to track down an explicit set of algebraically independent generators for $\mathbb D(\mathrm{SPos}_n(\R))$ in the literature. In \cite[\S~2]{Ga76} one finds a set of algebraically independent generators for the center $\mathcal Z$ of the universal enveloping algebra of $\mathfrak{sl}_n(\C)$, which is known to be isomorphic to the algebra of invariant differential operators on $\SL_n(\R)/\SO_n(\R)$ (see e.g. \cite[Prop.~II.5.32\,\&\,Exer.~II.D.3]{He84}). In \cite[\S~7]{Os08} the author follows a similar line giving the generators (without proof) as coefficients in a characteristic polynomial.

In this paper we show how to derive algebraically independent generators for $\mathbb D(\mathrm{SPos}_n(\R))$ from the Maass-Selberg operators as the $\mathrm{SPos}_n(\R)$-parts (constructed in Lemma~\ref{DiffOpProjection}) of the Maass-Selberg operators $\delta_2,\ldots,\delta_n$ (Theorem~\ref{main result}). The arguments given can be applied in more generality as we will explain in Section~\ref{sec:related results}. The key input of the paper comes from the master theses \cite{Br20,Ci16} written by the first two authors under the direction of the third author.

\section{Preliminaries on Invariant Differential Operators}

For a real reductive Lie group $G$ in the Harish-Chandra class \cite{GV88}  and a maximal compact subgroup $K \le G$ we denote by
$$\D_K(G)=  \{D \in \D(G) \mid \forall\, k \in K: \mathrm{Ad}(k)D=D \}$$
the $K$-invariant elements of $\D(G)$, i.e., the differential operators that are $G$-invariant from the left and $K$-invariant from the right. We also declare
$$\D(G)\mathfrak{k} = \braket{D\tilde{X} \mid D \in \D(G),X \in \mathfrak{k}}_{\C-\text{vector space}}$$
to be the left ideal generated by $\mathfrak{k}$ in $\D(G)$. Here $\tilde X$ is the left-invariant vector field on $G$ associated with $X$. 
\begin{Proposition}[{\cite[Thm.~II.4.6]{He84}}]\label{ReductiveDiffOps}
Let $\pi:G \to G/K$ be the canonical projection. Then the map
$$\mu:\D_K(G) \to \D(G/K), \; (\mu(D)f)(gK) = D(f\circ \pi)(g)  , \quad f \in C^\infty(G/K), g \in G$$
is an algebra epimorphism with kernel $\D_K(G) \cap \D(G)\mathfrak{k}$. Moreover, we have an isomorphism of algebras
$$\mu:\D_K(G)/(\D_K(G) \cap \D(G)\mathfrak{k}) \to \D(G/K).$$
\end{Proposition}

There exists a unique linear isomorphism, called the symmetrization of $\D(G)$,
$$\lambda_G = \lambda:\mathcal{S}(\mathfrak{g}_\C) \to \D(G)$$
satisfying $\lambda(X^m)=\tilde{X}^m$ for each $X \in \mathfrak{g}$, where $\mathcal{S}(\mathfrak{g}_\C)$ is the symmetric algebra of $\mathfrak{g}_\C$. If $X_1,\ldots,X_n$ is a basis of $\mathfrak{g}$, then  under the identification $\C[X_1,\ldots,X_n]\cong\mathcal{S}(\mathfrak{g}_\C)$, we obtain for each $P \in \mathcal{S}(\mathfrak{g}_\C)$
\begin{equation}\label{DefSymmetrizationDiffOp}
\lambda(P)f(g)=\left.P\left(\frac{\partial}{\partial t_1},\ldots,\frac{\partial}{\partial t_n}\right)\right|_{t=0} f(g\exp(t_1X_1+\ldots+t_nX_n)) , \; f\in C^\infty(G), g \in G, 
\end{equation}
where the suffix $|_{t=0}$ means the evaluation in $t=(t_1,\ldots,t_n)=0$ after differentiation (\cite[Thm.~II.4.3]{He84}).
Furthermore, for any $Y_1,\ldots,Y_r \in \mathfrak{g}$ we have
\begin{equation}\label{SymChara}
\lambda(Y_1\cdots Y_r)=\frac{1}{r!}\sum\limits_{\sigma \in \mathcal{S}_r} \tilde{Y}_{\sigma(1)}\cdots \tilde{Y}_{\sigma(r)},
\end{equation}
where $\mathcal{S}_r$ is the symmetric group of permutations on $r$ elements.

Let $\mathfrak{g}=\mathfrak{k} \oplus \mathfrak{p}$ be the Cartan decomposition and $I(\mathfrak p_\C)$ the space of complex $K$-invariant polynomials on $\mathfrak p_\C$.

\begin{Proposition}[{\cite[Thm.~II.4.9]{He84}}]\label{DiffOpLinearize}
 The map
$$\mu \circ \lambda_G : I(\mathfrak{p}_\C) \to \D(G/K), \quad P \mapsto D_P = \mu(\lambda_G(P))$$
is a linear isomorphism. Moreover, for each  $P_1,P_2 \in I(\mathfrak{p}_\C)$ there exists a $Q \in I(\mathfrak{p}_\C)$ with degree smaller than the sum of the degrees of $P_1,P_2$ and
$$D_{P_1P_2}=D_{P_1}D_{P_2} + D_Q.$$
\end{Proposition}

From this we can see that, if $P_1,\ldots, P_n$ are generators of $I(\mathfrak{p}_\C)$, then $D_{P_1},\ldots, D_{P_{n}}$ are generators of $\D(G/K)$.

\begin{Remark}\label{K-Invariance}
If $D \in \D(G/K)$ and $X_1,\ldots,X_r$ is a basis of $\mathfrak{p}$, then there exists a polynomial $P \in \mathcal{S}(\mathfrak{p}_\C)$, such that
$$Df(gH) = \left.P\left( \frac{\partial}{\partial t_1},\ldots, \frac{\partial}{\partial t_r}\right)\right|_{t=0} f(g\exp(t_1X_1+\ldots+t_rX_r) K), \quad f \in C^\infty(G/K),$$
which is uniquely determined, since $(t_1,\ldots,t_r) \mapsto  g\exp(t_1X_1+\ldots+t_rX_r)K$ is a local chart around $gK$. 
Thus, Proposition~\ref{DiffOpLinearize} implies that $P$ is automatically $K$-invariant, i.e., $P \in I(\mathfrak{p}_\C)$. 
\end{Remark}

\section{Maass-Selberg Operators and their $\mathrm{SPos}_n(\R)$-Radial Parts}\label{sec:radial parts}

Let $\mathrm{Sym}_n(\R)$ be the set of symmetric $n\times n$-matrices. Then the $k$-th \emph{Maass-Selberg operator} $\delta_k\in \D(\mathrm{Pos}_n(\R))$ is given by 
$$\delta_k f(g \mathrm{O}_n(\R)) = \left.\tr\left(\left( \frac{\partial}{\partial X}\right)^k\right)\right|_{X=0} f(g\exp(X)  \mathrm{O}_n(\R)) , \; f \in C^\infty(\GL_n(\R)/\mathrm{O}_n(\R)),$$
where $$\frac{\partial}{\partial X} = \begin{pmatrix}
\frac{\partial}{\partial x_{11}} & \hdots & \frac{1}{2}\frac{\partial}{\partial x_{1n}} \\
\vdots & \ddots & \vdots\\
\frac{1}{2}\frac{\partial}{\partial x_{1n}} & \hdots & \frac{\partial}{\partial x_{nn}} 
\end{pmatrix} \text{ and }X = \begin{pmatrix}
x_{11} & \hdots & x_{1n} \\
\vdots & \ddots & \vdots\\
x_{1n} & \hdots & x_{nn} 
\end{pmatrix} \in \mathrm{Sym}_n(\R).$$

To construct the $\mathrm{SPos}_n(\R)$-parts of the $\delta_k$ we consider the maps
$$i:\SL_n(\R)/\mathrm{SO}_n(\R) \to \GL_n(\R)/\mathrm{O}_n(\R) , \quad g\mathrm{SO}_n(\R) \mapsto g\mathrm{O}_n(\R)$$
and
$$p:\GL^+_n(\R)/\mathrm{SO}_n(\R) \to \SL_n(\R)/\mathrm{SO}_n(\R) , \quad g\mathrm{O}_n(\R) \mapsto \det(g)^{-\frac{1}{n}}g\mathrm{SO}_n(\R),$$
where $\GL_n^+(\R)= \{g \in \GL_n(\R) \mid \det(g)>0\}$. Note that 
$$\GL^+_n(\R)/\mathrm{SO}_n(\R)\to \GL_n(\R)/\mathrm{O}_n(\R),\quad g\mathrm{SO}_n(\R)\mapsto g\mathrm{O}_n(\R)$$
is a diffeomorphism. We use it to identify the two spaces. For $g \in \SL_n(\R)$ let $\tau_g$ denote the translation by $g$ on $\GL_n(\R)/\mathrm{O}_n(\R)$ and $\SL_n(\R)/\mathrm{SO}_n(\R)$.

\begin{Proposition}\label{TransProjInj}
The maps $i$ and $p$ are $\SL_n(\R)$-equivariant, i.e. 
$$\forall g\in \SL_n(\R):\quad  p \circ \tau_g=\tau_g \circ p \text{ and } i \circ \tau_g=\tau_g \circ i.$$
\end{Proposition}

\begin{Proof}
Let $h\mathrm{O}_n(\R) \in \GL_n(\R)/\mathrm{O}_n(\R)$ and $g \in \SL_n(\R)$. We can assume that $\det(h)>0$, so
$$p(g.h\mathrm{O}_n(\R))=\det(gh)^{-\frac{1}{n}}gh\mathrm{SO}_n(\R)=g.\det(h)^{-\frac{1}{n}}h\mathrm{SO}_n(\R) = g.p(h\mathrm{O}_n(\R)).$$
For $h\mathrm{SO}_n(\R) \in \SL_n(\R)/\mathrm{SO}_n(\R)$ we have
$$i(g.h\mathrm{SO}_n(\R))=gh\mathrm{O}_n(\R)=g.h\mathrm{O}_n(\R)=g.i(h\mathrm{SO}_n(\R)).$$
\end{Proof}

\begin{Lemma}\label{DiffOpProjection}
The map
$P:\D(\GL_n(\R)/\mathrm{O}_n(\R)) \to \D(\SL_n(\R)/\mathrm{SO}_n(\R))$
defined by
$$ P(D)f = D(f\circ p) \circ i$$
for $f\in C^\infty(\SL_n(\R)/\mathrm{SO}_n(\R))$ is a morphism of algebras. Furthermore, it  satisfies
\begin{equation}\label{FactorizingDiffOp}
(P(D)f) \circ p = D(f \circ p)
\end{equation}
for $D\in \D(\GL_n(\R)/\mathrm{O}_n(\R))$ and $f \in C^\infty(\SL_n(\R)/\mathrm{SO}_n(\R))$.
\end{Lemma}

We call $P(D)$ the \emph{$\mathrm{SPos}_n(\R)$-radial part} of $D$.

\begin{Proof}
First we check that the image of $P$ consists of invariant differential operators on $\SL_n(\R)/\mathrm{SO}_n(\R)$. If $f \in C^\infty(\SL_n(\R)/\mathrm{SO}_n(\R))$, then  $f \circ p$ is smooth since $p$ is. As $D$ is a differential operator and $i$ is smooth, we have $P(D)f=D(f\circ p) \circ i \in C^\infty(\SL_n(\R)/\mathrm{SO}_n(\R))$. The linearity of $P(D)$ is an easy consequence of the linearity of $D$. Moreover, $D$ satisfies $\mathrm{supp}(D(f\circ p)) \subseteq \mathrm{supp}(f\circ p)$. Precomposing by  $i$ we find
$$\mathrm{supp}(P(D)f)=\mathrm{supp}(D(f\circ p) \circ i) \subseteq \mathrm{supp}(f\circ p \circ i) = \mathrm{supp}(f),$$
so that Peetre's Theorem shows that $P(D)$ is a differential operator on $\SL_n(\R)/\mathrm{SO}_n(\R)$. From the invariance of $D \in \D(\GL_n(\R)/\mathrm{O}_n(\R))$ and Proposition~\ref{TransProjInj} we obtain
$$P(D)(f\circ \tau_g) = D(f\circ \tau_g \circ p) \circ i = D(f\circ p) \circ i \circ \tau_g=P(D)f \circ \tau_g$$
for $g\in\SL_n(\R)$, so $P(D)\in\D(\SL_n(\R)/\mathrm{SO}_n(\R))$.

The linearity of $P$ is clear. Next we prove equation \eqref{FactorizingDiffOp}.  Let $g\mathrm{O}_n(\R) \in \GL_n(\R)/\mathrm{O}_n(\R)$ be arbitrary and assume that $\det(g)>0$. Define $h = \det(g)^{-\frac{1}{n}}\1_n$. Then by the invariance of $D \in \D(\GL_n(\R))/\mathrm{O}_n(\R))$ and Proposition~\ref{TransProjInj} we obtain
\begin{align*}
(P(D)f \circ p)(g\mathrm{O}_n(\R)) &= (D(f\circ p) \circ i)(\det(g)^{-\frac{1}{n}}g\mathrm{SO}_n(\R)) \\
&= D(f \circ p)(\det(g)^{-\frac{1}{n}}g\mathrm{O}_n(\R)) \\
&= (D(f\circ p)\circ \tau_h)(g\mathrm{O_n}(\R)) \\
&= D(f\circ p \circ \tau_h)(g\mathrm{O}_n(\R)).
\end{align*}
Since $\tau_h$ only scales the determinant we have $p\circ \tau_h=p$, this implies \eqref{FactorizingDiffOp}. 

Let $D_1,D_2 \in \D(\GL_n(\R)/\mathrm{SL}_n(\R))$. Together with equation \eqref{FactorizingDiffOp}, we conclude
$$P(D_1D_2)f = (D_1D_2(f\circ p)) \circ i = D_1(P(D_2)f\circ p) \circ i = P(D_1)P(D_2)f.$$
As $P(\mathrm{id}_{\D(\GL_n(\R)/\mathrm{O}_n(\R)}))=\mathrm{id}_{\D(\SL_n(\R)/\mathrm{SO}_n(\R))}$ is obvious, this concludes the proof.
\end{Proof}

By Proposition \ref{DiffOpLinearize} for each differential operator $D \in \D(\GL_n(\R)/\mathrm{O}_n(\R))$ there exists  a polynomial $Q_D$ on $\mathrm{Sym}_n(\R)$, such that
$$Df(g\mathrm{O}_n(\R)) = \left. Q_D\left( \frac{\partial}{\partial X}\right)\right|_{X=0}f(g\exp(X)\mathrm{O}_n(\R)),$$
where $X \in \mathrm{Sym}_n(\R)$ and $f \in C^\infty(\GL_n(\R)/\mathrm{O}_n(\R))$. This leads to another representation of the morphism $P$.

\begin{Proposition}\label{ProjectionImage}
For each $f \in C^\infty(\SL_n(\R)/\mathrm{SO}_n(\R))$
$$P(D)f(g\mathrm{SO}_n(\R)) = \left. Q_D\left( \frac{\partial}{\partial X}\right)\right|_{X=0}f\left(g\exp\left(X-\frac{1}{n}\tr(X)\1_n\right)\mathrm{SO}_n(\R)\right).$$
\end{Proposition}

\begin{Proof}
\begin{align*}
P(D)f(g\mathrm{SO}_n(\R)) 
&= D(f\circ p)(g\mathrm{O}_n(\R))) \\
&= \left.Q_D\left( \frac{\partial}{\partial X}\right)\right|_{X=0} f(p(g\exp(X)\mathrm{O}_n(\R))) 
\\
&=\left.Q_D\left( \frac{\partial}{\partial X}\right)\right|_{X=0} f(\det(\exp(X))^{-\sfrac{1}{n}}g\exp(X)\mathrm{SO}_n(\R))  \\
&= \left.Q_D\left( \frac{\partial}{\partial X}\right)\right|_{X=0} f\left(g\exp\left(-\frac{1}{n}\tr(X)\1_n\right)\exp(X)\mathrm{O}_n(\R)\right) \\
&= \left.Q_D\left( \frac{\partial}{\partial X}\right)\right|_{X=0} f\left(g\exp\left(X-\frac{1}{n}\tr(X)\1_n\right) \cdot \mathrm{SO}_n(\R)\right). 
\end{align*}
\end{Proof}

In fact, the next lemma shows that $P$ is surjective. Thus, the images $P(\delta_1),\ldots, P(\delta_n)$ of the Maas-Selberg operators are generators of $\D(\mathrm{SPos}_n(\R))$ and we will show that $P(\delta_1)=0$, so that it remains to show that $P(\delta_2),\ldots,P(\delta_n)$ are algebraically independent.

\begin{Lemma}\label{Surjectivity}
The morphism $P$ is surjective and $P(\delta_1)=0$.
\end{Lemma}

\begin{Proof}
Let $D \in \D(\mathrm{Pos}_n(\R))$ be given by a unique polynomial $Q$ on $\mathrm{Sym}_n(\R)$, i.e.
$$Df(g\mathrm{O}_n(\R))=\left.Q\left(\frac{\partial}{\partial X}\right)\right|_{X=0} f(g\exp(X)\mathrm{O}_n(\R))$$
holds for every $f \in C^\infty(\mathrm{GL}_n(\R)/\mathrm{O}_n(\R))$. Denote by $\braket{X,Y}= \tr(XY)$ the inner product on $\mathrm{Sym}_n(\R)$.
Then it is easy to verify that for $A \in \mathrm{Sym}_n(\R)$ we have 
\begin{equation}\label{DiffSymExp}
\frac{\partial}{\partial x_{ii}} e^{\braket{X,A}} = a_{ii}e^{\braket{X,A}}, \quad \frac{1}{2}\frac{\partial}{\partial x_{ij}} e^{\braket{X,A}}= a_{ij}e^{\braket{X,A}}.
\end{equation}
Let $Q_1$ be the unique polynomial on $\mathrm{SSym}_n(\R)$ so that for each function $f \in C^\infty(\SL_n(\R)/\mathrm{SO}_n(\R))$ we have
$$P(D)f(g\mathrm{SO}_n(\R))=\left.Q_1\left(\frac{\partial}{\partial Y}\right)\right|_{Y=0} f(g\mathrm{exp}(Y)\mathrm{SO}_n(\R)).$$
We associate to a matrix $A \in \mathrm{SSym}_n(\R)$ the smooth function 
$$f_A:\GL_n(\R)/\mathrm{O}_n(\R) \to \R, \quad \exp(Y)\mathrm{SO}_n(\R) \mapsto e^{\braket{Y,A}} \text{ with } Y \in \mathrm{SSym}_n(\R).$$ 
We can observe with Proposition \ref{ProjectionImage} and equation \eqref{DiffSymExp} that
\begin{align*}
Q(A) &= \left.Q\left( \frac{\partial}{\partial X}\right)\right|_{X=0} e^{\braket{X,A}} = \left.Q\left( \frac{\partial}{\partial X}\right)\right|_{X=0} e^{\tr\left(\left(X-\frac{\tr(X)}{n}\1_n\right)A\right)} \\
&= \left.Q\left( \frac{\partial}{\partial X}\right)\right|_{X=0} f_A\left(\exp\left(X-\frac{\tr(X)}{n}\1_n\right)\mathrm{SO}_n(\R)\right) \\
&= P(D)f_A(\1_n\mathrm{SO}_n(\R)) \\
&= \left.Q_1\left( \frac{\partial}{\partial Y}\right)\right|_{Y=0} f_A\left(\exp(Y)\mathrm{SO}_n(\R)\right) \\
&= \left.Q_1\left( \frac{\partial}{\partial Y}\right)\right|_{Y=0} e^{\braket{Y,A}} = Q_1(A).
\end{align*}
Thus $Q_1$ is the restriction of $Q$ to $\mathrm{SSym}_n(\R)$, so $P$ is surjective. Moreover, for $D=\delta_1$ we have $Q=\tr$ and therefore $Q_1=0$, hence $P(\delta_1)=0$.
\end{Proof}

To summarize, we have shown the following theorem.

\begin{Theorem}\label{SLnDiffOpsTrace}
Let $\delta_1,\ldots, \delta_n$ be the Maass-Selberg operators. Then  $P(\delta_1)=0$ and 
$$
P(\delta_k)f(g\mathrm{SO}_n(\R)) = \left.\tr\left( \left(\frac{\partial}{\partial X}\right)^k\right)\right|_{X=0} f\left(g\exp\left(X-\frac{1}{n}\tr(X)\1_n\right) \cdot \mathrm{SO}_n(\R)\right),
$$
where 
$$\frac{\partial}{\partial X} = \begin{pmatrix}
\frac{\partial}{\partial x_{11}} & \hdots & \frac{1}{2}\frac{\partial}{\partial x_{1n}} \\
\vdots & \ddots & \vdots\\
\frac{1}{2}\frac{\partial}{\partial x_{1n}} & \hdots & \frac{\partial}{\partial x_{nn}} 
\end{pmatrix} \text{ and }X=\begin{pmatrix}
x_{11} & \hdots & x_{1n} \\
\vdots & \ddots & \vdots\\
x_{1n} & \hdots & x_{nn} 
\end{pmatrix} \in \mathrm{Sym}_n(\R)$$
are generators for the algebra $\D(\SL_n(\R)/\mathrm{SO}_n(\R))$.
\end{Theorem}

Using the proof of Lemma \ref{Surjectivity} we can express the $P(\delta_k)$ in terms of a local chart for $\mathrm{SPos}_n(\R)=\SL_n(\R)/\mathrm{SO}_n(\R)$.

\begin{Corollary}
For $f \in C^\infty(\SL_n(\R)/\mathrm{SO}_n(\R))$ and $Y \in \mathrm{SSym}_{n}(\R)$ we have
\begin{equation}\label{MaassSelberSLnChart}P(\delta_k)f(g\mathrm{SO}_n(\R)) = \left.\tr\left( \left(\frac{\partial}{\partial Y}\right)^k\right)\right|_{Y=0} f(g\exp(Y)\cdot \mathrm{SO}_n(\R)).
\end{equation}
\end{Corollary}

It is possible to prove the surjectivity of $P$ by giving an explicit algebra splitting. In fact, for $g \mathrm{SO}_n(\R) \in \GL_n^+(\R)/\mathrm{SO}_n(\R)$ one can define the map
$$\varphi_g:\GL_n^+(\R)/\mathrm{SO}_n(\R) \to \GL_n^+(\R)/\mathrm{SO}_n(\R), \quad h\mathrm{SO}_n(\R) \to \det(g)^{\frac{1}{n}}h\mathrm{SO}_n(\R),$$
which only depends on the coset $g\mathrm{O}_n(\R)$. Using $\varphi_g$ and the identification $\GL_n(\R)/\mathrm{O}_n(\R)=\GL_n^+(\R)/\mathrm{SO}_n(\R)$ one obtains an algebra morphism
$I:\D(\mathrm{SPos}_n(\R)) \to \D(\GL_n(\R)/\mathrm{O}_n(\R))$ via
$$I(D)f(g\mathrm{O}_n(\R)) = (D(f\circ \varphi_g \circ i) \circ p)(g\mathrm{O}_n(\R))$$
for $f \in C^\infty(\GL_n(\R)/\mathrm{O}_n(\R))$ and $D \in \D(\SL_n(\R)/\mathrm{SO}_n(\R))$. This morphism satisfies
$$I(D)f\circ i = D(f\circ i), \quad P \circ I = \mathrm{id}_{\D(\SL_n(\R)/\mathrm{SO}_n(\R))}.$$
For more details we refer to \cite[Lemma 4.5, Theorem 4.6]{Br20}.

\section{Algebraic Independence}

In this section we prove the algebraic independence of $P(\delta_2),\ldots, P(\delta_n)$. We start with a general lemma on polynomial algebras.

\begin{Lemma}\label{CommutativeGeneratorOfAlgebra}
Let $k$ be a field and $A$ an (associative) commutative unital algebra over $k$. Assume that $A$ contains algebraically independent generators $x_1,\ldots, x_n \in A$. Then
\begin{enumerate}
\item[(i)] If $y_1,\ldots,y_m \in A$ generate $A$, then $m \ge n$.
\item[(ii)] If $y_1,\ldots, y_n \in A$ generate $A$, then $y_1,\ldots,y_n$ are algebraically independent.
\item[(iii)] The number $n$ is independent of the choice of algebraically independent generators.
\end{enumerate}
\end{Lemma}

\begin{Proof}
Since $x_1,\ldots,x_n \in A$ are algebraically independent generators of $A$, we have an algebra isomorphism
$$k[X_1,\ldots,X_n] \to A, \quad X_i \mapsto x_i.$$
So we may assume that $A=k[X_1,\ldots,X_n]$ is the ring of polynomials in $n$-indeterminants and $x_i=X_i$.
\begin{enumerate}
\item[(i)] We define a morphism of algebras by
$$k[Y_1,\ldots,Y_m] \to k[X_1,\ldots,X_n], \quad Y_i \mapsto y_i.$$
The morphism $\varphi$ is surjective, because $y_1,\ldots,y_m$ are generators of $A=k[X_1,\ldots,X_n]$. If $\dim$ denotes the Krull-dimension (cf. \cite[Def.~2.5.3]{Bo19}), then the surjectivity of $\varphi$ implies $\dim(k[Y_1,\ldots,Y_m]) \ge \dim(k[X_1,\ldots,X_n])$ as preimages of prime ideals are prime (cf. \cite[Prop.~2.5.5]{Bo19}). Now \cite[Corollary 2.25.2]{Bo19} implies 
$\dim(k[Y_1,\ldots,Y_m])=m$ and $\dim(k[X_1,\ldots,X_n]) = n$, and hence the claim.

\item[(ii)] If $y_1,\ldots,y_n$ are generators of $A=k[X_1,\ldots,X_n]$, then it suffices to show that the map $\varphi$ from part (i) is injective, because the $Y_1,\ldots,Y_n$ are algebraically independent.

Assume that $0 \neq \ker(\varphi)$. Then we can find a polynomial $0 \neq f \in \ker(\varphi)$. Decompose $f=f_1\cdots f_r$ into irreducible polynomials $f_i$. Then 
$$0=\varphi(f)=\varphi(f_1)\cdots \varphi(f_r)$$ leads to an irreducible polynomial $f_i \in \ker(\varphi)$, since $k[X_1,\ldots,X_n]$ has no zero divisors. 
Let $I\subseteq \ker(\varphi)$ be the ideal generated by $f_i$. Moreover, $I$ is a prime ideal which does not contain any proper prime ideal except $\{0\}$. This means that the height $\mathrm{ht}(I)$ of $I$ is $1$ (cf. \cite[Def.~2.5.4\,\&\,Prop.~2.25.3]{Bo19}), so that \cite[Lemma 2.25.7]{Bo19} implies
$\dim(k[Y_1,\ldots,Y_n]/I)=n-1$. 
The morphism $\varphi$ factors through $I$ to a  surjective morphism
$\tilde{\varphi}:k[Y_1,\ldots,Y_n]/I \to k[X_1,\ldots,X_n]$. Again, we conclude 
$\dim(k[Y_1,\ldots,Y_n]/I)\ge \dim(k[X_1,\ldots,X_n])$ and obtain a contradiction proving the claim.
\item[(iii)] This follows applying parts (i) and (ii) to two different algebraically independent generators $x_1,\ldots,x_n$ and $x'_1,\ldots,x_m'$ of $A$.
\end{enumerate} 
\end{Proof}

\begin{Proposition}\label{SLnDiffOpsTrace2}
Let $\delta_1,\ldots, \delta_n$ be the Maass-Selberg operators. Then the operators $P(\delta_2),\ldots,P(\delta_2)$ are algebraically independent.
\end{Proposition}

\begin{Proof}
By the Harish-Chandra Isomorphism and Chevalley's Theorem about invariants of finite reflection groups,  we know that $\D(\SL_n(\R)/\mathrm{SO}_n(\R))$ is generated by $n-1$ algebraically independent generators. Since $\delta_1,\ldots, \delta_n$ are generators of $\D(\GL_n(\R)/\mathrm{O}_n(\R))$ and $P:\D(\GL_n(\R)/\mathrm{O}_n(\R)) \to\D(\SL_n(\R)/\mathrm{SO}_n(\R))$ is an epimorphism of algebras and $P(\delta_1)=0$ by Theorem~\ref{SLnDiffOpsTrace}, the operators $P(\delta_2),\ldots, P(\delta_n)$ generate $\D(\SL_n(\R)/\mathrm{SO}_n(\R))$. Finally, Lemma~ \ref{CommutativeGeneratorOfAlgebra} implies the algebraic independence.
\end{Proof}

Combining Proposition~\ref{SLnDiffOpsTrace2} with the results of Section~\ref{sec:radial parts} we obtain

\begin{Theorem}\label{main result}
The operators $P(\delta_2),\ldots,P(\delta_2)$ given by \eqref{MaassSelberSLnChart} form an algebraically independent set of generators for $\D(\mathrm{SPos}_n(\R))$.
\end{Theorem}

\section{Related Results}
\label{sec:related results}

The arguments given in this paper also work for different sets of algebraically independent generators for  $\D(\mathrm{Pos}_n(\R))$. 

\begin{Remark}[{\cite[Thm.~3.7\&Prop.~4.13]{Br20}}]
Let $X=(x_{ij})_{i,j} \in \mathfrak{p}$ be a symmetric matrix with eigenvalues $\lambda_1,\ldots,\lambda_n$. Assume that $Y_1,\ldots, Y_n$ are the elementary symmetric polynomials in $n$ indeterminants. Then one can see that the characteristic polynomial of $X$ can be written as
\begin{equation}\label{CharaPoly}
\det(t\1_n-X)=\prod\limits_{i=1}^n (t-\lambda_i) = t^n+\sum\limits_{k=1}^n (-1)^kY_k(\lambda_1,\ldots,\lambda_n)t^{n-k}.
\end{equation}
Let $F_k(X)$ be the sum over the $k\times k$ principal minors of $X$, i.e.,
$$F_k(X)=\sum\limits_{1\le i_1<\ldots < i_k\le n} \det\left((x_{i_{j},i_{\ell}})_{j,\ell=1,\ldots,k}\right)=\sum\limits_{1\le i_1<\ldots < i_k\le n} \sum\limits_{\sigma \in \mathcal{S}_k} \mathrm{sgn}(\sigma)\prod\limits_{j=1}^k x_{i_ji_{\sigma(j)}},$$
where $\mathrm{sgn}(\sigma)$ is the sign of $\sigma$. Now we also obtain (cf. \cite[Formula 1.2.13]{HJ13}), that
\begin{equation}\label{MinorSum}
\det(t\1_n-X)= t^n+\sum\limits_{k=1}^n (-1)^kF_k(X)t^{n-k}.
\end{equation}
Thus, we can see that for each orthogonal matrix $k \in K$ we have
\begin{equation}\label{KInvariance}
F_j(kXk^{-1})=F_j(X) , \quad j=1,\ldots,n,
\end{equation}
since the characteristic polynomial is invariant under conjugations. Finally, combining \eqref{CharaPoly} and \eqref{MinorSum},
\begin{equation}\label{SymPolyAsMinor}
F_j(X)=Y_j(\lambda_1,\ldots,\lambda_n)  , \quad j=1,\ldots,n.
\end{equation}
The operators $\eta_1,\ldots, \eta_n$ given by
$$\eta_kf(g\mathrm{O}_n(\R)) = \left.F_k\left( \frac{\partial}{\partial X}\right)\right|_{X=0} f(g\exp(X)  \mathrm{O}_n(\R)) , \quad f \in C^\infty(\GL_n(\R)/\mathrm{O}_n(\R))$$
form an algebraically independent set of generators for $\D(\mathrm{Pos}_n(\R)$. Then
$$P(\eta_k)f(g\mathrm{SO}_n(\R)) = \left.F_k \left(\frac{\partial}{\partial Y}\right)\right|_{Y=0} f(g\exp(Y)\cdot \mathrm{SO}_n(\R))$$
holds for $f \in C^\infty(\SL_n(\R)/\mathrm{SO}_n(\R))$ and $Y \in \mathrm{SSym}_{n}(\R)$.
Moreover, $P(\eta_1)=0$ and $P(\eta_2),\ldots,P(\eta_n)$ form an algebraically independent set of generators for $\D(\mathrm{SPos}_n(\R))$.
\end{Remark} 

Similar results  hold for the symmetric spaces $\mathrm{Pos}_n(\C)=\GL_n(\C)/\mathrm{U}_n(\C)$ of positive definite hermitian $n\times n$-matrices and $\mathrm{SPos}_n(\C)=\SL_n(\C)/\mathrm{SU}_n(\C)$ of determinant $1$ elements in $\mathrm{Pos}_n(\C)$ (see \cite[Theorems ~3.8, 3.9, 4.17]{Br20}).

Note that $\mathrm{Pos}_n(\R)$ and $\mathrm{Pos}_n(\C)$ are examples for symmetric cones (see \cite{FK94}). Nomura gave algebraically independent generators $D_1,\ldots, D_r$ for the algebra of invariant differential operators on any symmetric cone $\Omega$ in \cite{No89}. Using the Jordan theoretic determinant one obtains the symmetric space $S\Omega$ associated with the derived groups $\mathrm{Str}(V)'$ of the structure group $\mathrm{Str}(V)$ of the euclidean Jordan algebra $V$ defined by $\Omega$ (see \cite{FK94} for the relevant definitions and constructions). The arguments presented in this paper can also be adapted to this situation to yield an algebra epimorphism $P:\D(\Omega)\to \D(S\Omega)$  via $\mathrm{Str}(V)'$-radial parts such that $\C D_1=\ker P$. One obtains the following theorem.

\begin{Theorem}
Let $\Omega$ be a symmetric cone and $D_1,\ldots,D_r$ the set of Nomura's generators for $\D(\Omega)$. Then $P(D_2),\ldots,P(D_r)$ are algebraically independent generators for $\D(S\Omega)$.
\end{Theorem}   

The classification of symmetric cones (\cite[p.~97]{FK94}) shows that this yields algebraically independent generators for the algebras of invariant differential operators for the following Riemannian symmetric spaces: $\SL_n(\mathbb R)/\SO_n(\R)$, $\SL_n(\mathbb C)/\mathrm{SU}_n(\C)$,  $\SL_n(\mathbb H)/\mathrm{SU}_n(\H)$, $\mathrm O_{1,n-1}(\R)/\mathrm O_{n-1}(\R)$, $E_{6(-26)}/F_4$, where $\mathbb H$ denotes the quaternions and $E_{6(-26)}/F_4$ is the space of positive definite $3\times 3$-matrices of determinant $1$ with entries in the octonians.

\bigskip
\noindent
{\it Acknowledgment:} We thank the anonymous referee for suggesting the short argument leading to Theorem~\ref{SLnDiffOpsTrace}.

\end{document}